\documentclass[12pt]{article}
\everymath{\displaystyle}
\usepackage{graphicx}
\usepackage{pst-all}
\usepackage{color}
\usepackage{amssymb}
\usepackage{amsmath}
\newtheorem{prethm}{{\bf Theorem}}

\newenvironment{thm}{\begin{prethm}{\hspace{-0.5
               em}{\bf .}}}{\end{prethm}}
\newtheorem{prelemma}{{\bf Lemma}}

\newenvironment{lemma}{\begin{prelemma}{\hspace{-0.5
               em}{\bf .}}}{\end{prelemma}}
\newtheorem{preex}{{\bf Example}}

\newtheorem{preprop}{{\bf Proposition}}

\newenvironment{prop}{\begin{preprop}{\hspace{-0.5em}{\bf .}}}{\end{preprop}}
\newtheorem{precor}{{\bf Corollary}}

\newtheorem{preremark}{{\bf Remark}}

\newtheorem{preprob}{{\bf Problem}}

\newtheorem{predefin}{{\bf Definition}}

\newenvironment{defin}{\begin{predefin}{\hspace{-0.5
               em}{\bf .}}}{\end{predefin}}
\newtheorem{preconj}{{\bf Conjecture}}

\newtheorem{preprobb}{{\bf Problem}}

\newtheorem{prelem}{{\bf Theorem}}

\newenvironment{proof}{{\bf Proof.}\rm }{\hfill{$\Box$}}

\newtheorem{presolution}{{\bf Solution.}}

\def\newpic#1{}

\title{\Large\bf More results on greedy defining sets}
\author{
{\large\bf Manouchehr Zaker}
\vspace{3mm}\\
    {Institute for Advanced Studies in Basic Sciences, Zanjan, Iran}\\
    mzaker@iasbs.ac.ir}
\date{}
\begin{document}
\maketitle

\begin{abstract}
\noindent The greedy defining sets of graphs were appeared first
time in [M. Zaker, Greedy defining sets of graphs, Australas. J.
Combin, 2001]. We show that to determine the greedy defining number
of bipartite graphs is an NP-complete problem. This result answers
affirmatively the problem mentioned in the previous paper. It is
also shown that this number for forests can be determined in
polynomial time. Then we present a method for obtaining greedy
defining sets in Latin squares and using this method, show that any
$n\times n$ Latin square has a GDS of size at most $n^2-(n\log
n)/4$. Finally we present an application of greedy defining sets in
designing practical secret sharing schemes.
\end{abstract}


\section{Introduction}

\noindent Let $G$ be a simple graph whose vertices are ordered by an
order $\sigma$ as $v_1, \ldots, v_n$. The first-fit (greedy)
coloring of $G$ with respect to $\sigma$ starts with $v_1$ and
assigns color 1 to $v_1$ and then goes to the next vertex. It colors
$v_i$ by the first available color which is not appeared in the
neighborhood of $v_i$. If the algorithm finishes coloring of $G$ by
$\chi(G)$ colors then we say that it succeeds. But this is not the
case in general. If we want the greedy algorithm to succeed, then we
need to pre-color some of the vertices in $G$ before the algorithm
is invoked. So we define a {\it greedy defining set} (GDS) to be a
subset of vertices in $G$ together with a pre-coloring of $S$, that
will cause the greedy algorithm to successfully color the whole
graph $G$ with $\chi(G)$ colors. It is understood that the algorithm
skips over the vertices that are part of the defining set. Greedy
defining sets of graphs were first defined and studied by the author
in \cite{Z1}. This concept have also been studied for Latin squares
in \cite{Z2,Z3} and recently in \cite{V}. In the sequel follow the
formal definitions.

\begin{defin}
For a graph $G$ and an order $\sigma$ on $V(G)$, a {\it greedy
defining set} is a subset $S$ of $V(G)$ with an assignment of colors
to vertices in $S$, such that the pre-coloring can be extended to a
$\chi(G)$- coloring of $G$ by the greedy coloring of $(G,\sigma)$
and fixing the colors of $S$. The {\it greedy defining number} of
$G$ is the size of a greedy defining set which has minimum
cardinality, and is denoted by ${\rm GDN}(G,\sigma)$. A greedy
defining set for a $\chi(G)$-coloring $C$ of $G$ is a greedy
defining set of $G$ which results in $C$. The size of a greedy
defining set of $C$ with the smallest cardinality is denoted by
${\rm GDN}(G,\sigma,C)$.
\end{defin}

\noindent Let an ordered graph $(G,\sigma)$ and a proper vertex
coloring $C$ of $G$ using $\chi(G)$ colors be given. Let $i$ and $j$
with $1\leq i<j\leq \chi(G)$ be two arbitrary and fixed colors. Let
a vertex say $v$ of color $j$ be such that all of its neighbors with
color $i$ (this may be an empty set) are higher than $v$. Then $v$
together with these neighbors form a subset which we call a {\it
descent}. It was proved in \cite{Z1} that a subset $S$ of vertices
is a greedy defining set for the triple $(G,\sigma,C)$ if and only
if $S$ intersects any descent of $G$ or equivalently $S$ is a
transversal for the set of all descents.

\noindent Consider the Cartesian product $K_n\Box K_n$ and the
lexicographic order of its vertices. Namely $(i,j)<(i',j')$ if and
only if either $i<i'$ or $i=i'$ and $j<j'$. Since any $n\times n$
Latin square is equivalent to a proper $n$-coloring of the Cartesian
product $K_n\Box K_n$ then we can define greedy defining set and
number of Latin squares. In any Latin square we denote any cell in
row $i$, column $j$ with entry $x$ by $(i,j;x)$. Now the concept of
descent in the context of Latin squares is stated as follows. Given
a Latin square $L$, a set consisting of three cells $(i,j;y)$,
$(r,j;y)$ and $(i,k;x)$ where $i<r$, $j<k$ and $x<y$, is called a
descent. The following theorem proved in \cite{Z2,Z3} is in fact a
consequence of the theorem concerning GDS and transversal of
descents which was mentioned in the previous paragraph.

\begin{thm}
A subset $D$ of entries in a Latin square $L$ is greedy defining set
if and only if $D$ intersects any descent of $L$.\label{descent}
\end{thm}

\section{Greedy defining number of graphs}

\noindent In \cite{Z1} the computational complexity of determining
${\rm GDN}(G,\sigma,C)$ has been studied.

\begin{thm}{\rm (\cite{Z1})}
\noindent Given a triple $(G,\sigma,C)$ and an integer $k$. It is an
NP-complete problem to decide ${\rm GDN}(G,\sigma,C)\leq
k$.\label{2}
\end{thm}

\noindent Throughout the paper by the vertex cover problem we mean
the following decision problem. Given a simple graph $F$ and an
integer $k$, whether $F$ contain a vertex cover of at most $k$
vertices? Recall that a vertex cover is a subset $K$ of vertices
such that any edge is incident with a vertex of $K$. This problem is
a well-known NP-hard problem. In \cite{Z1} the vertex cover problem
was used to prove Theorem \ref{2}. But because there exists a flaw
in its proof, in the sequel we first fix the proof by slight
modification of it and then discuss the open question posed in
\cite{Z1}.

\noindent {\bf Proof of Theorem 2.} It is enough to reduce the
vertex cover problem to our problem. Let $(F,k)$ be an instance of
the vertex cover problem where $F$ has order $n$. We first color
arbitrarily the vertices of $F$ by $n$ distinct colors. Denote the
color of a vertex $v\in F$ by $c(v)$. Now we consider the complete
graph $K_n$ (vertex disjoint from $F$) on vertex set $\{1,2, \ldots,
n\}$. We order a vertex $i$ in $K_n$ by the very $i$ and a vertex
$v\in F$ by $2n-j+1$ if $c(v)=j$. For any $i$ and $j$ with $i<j$, we
put an edge between a vertex $v$ of $F$ of color $j=c(v)$ and a
vertex $i$ from $K_n$ if and only if $v$ is not adjacent to the
vertex of color $i$ in $F$. Let the color of a vertex $i\in K_n$ be
$i$. Denote the resulting ordered graph $(G,\sigma)$ and the proper
coloring of $G$ by $C$. It is easily checked that no descent in
$(G,\sigma,C)$ consists of only a single vertex. Since the colors of
$F$ are all distinct then a descent can only have two vertices and
we note that any edge in $F$ forms in fact a descent and these are
the only descents of $G$. We conclude that a transversal for the set
of descents in $G$ is a vertex cover for $F$ and vise versa. This
completes the proof.\\

\noindent It was asked in \cite{Z1} that given an ordered graph
$(G,\sigma)$, whether to determine $GDN(G,\sigma)$ is an NP-complete
problem? This problem is in fact the uncolored version of Theorem
\ref{2} where no coloring of graph is given in the input.

\noindent In the following we answer this problem affirmatively. We
begin with the following lemma.

\begin{lemma}\label{lem}
Let $G$ be a connected bipartite graph with a bipartition $(X,Y)$
whose vertices are colored properly by 1 and 2. Then there exists a
minimum greedy defining set for $G$ consisting of only vertices
colored 1.
\end{lemma}

\noindent \begin{proof} Let $S$ be a minimum GDS which contains the
minimum number of vertices colored 2. Consider a vertex $v\in S$ of
color 2. Since $S$ is a minimum GDS, $S\setminus \{v\}$ is not a GDS
and so all neighbors of $v$ with color 1 appear after $v$ in the
ordering of $G$. Now it suffices to delete $v$ from $S$ and add any
neighbor of it to $S$. The new member has color 1 because there are
only two colors in the graph. The resulting set is still a GDS and
this contradicts with our choice of $S$. Therefore there exists a
minimum GDS containing no vertex of color 2.
\end{proof}

\begin{thm}
Given an ordered connected bipartite graph $G$ and a positive
integer $k$. It is NP-complete to decide whether $GDN(G)\leq k$.
\end{thm}

\noindent \begin{proof} We transform an instance $(F,k)$ of the
vertex cover problem to an instance of our problem where $F$ is a
connected graph. Let $V(F)$ and $E(F)$ be the vertex and edge set of
$F$, respectively. Assume that $V_1$ and $V_2$ are two disjoint
copies of $V(F)$. Namely any vertex of $F$ has two distinct copies
in $V_1$ and $V_2$. Similarly let $E_1$ and $E_2$ be two disjoint
copies of $E(F)$. Let $G$ be the bipartite graph consisting of the
bipartite sets $X= V_1\cup E_1$ and $Y= V_2\cup E_2$, where $v\in
V_1\subseteq X$ is adjacent to $e\in E_2 \subseteq Y$ if $e$ (as an
edge of $F$) is incident to $v$ in $F$. Also a vertex $v\in V_2$ is
adjacent to $e\in E_1$ if $e$ is incident to $v$ in $F$. The only
extra edges of $G$ are of the form $vv'$ where $v$ is an arbitrary
vertex in $V_1$ and $v'$ its copy in $V_2$. We consider any ordering
$\sigma$ of $V(G)$ in which $E_2<E_1<V_2<V_1$, where for any two
sets $A$ and $B$ by $A<B$ we mean any element of $A$ has lower order
than any element of $B$.

\noindent The bipartite graph $G$ is connected since $F$ is so.
Therefore $G$ has only two proper colorings with two colors. To
determine $GDN(G)$ it is enough to determine the minimum greedy
defining number of these two colorings of $G$. Consider an arbitrary
coloring $C$ of $G$ in which the part $X$ is colored 1. According to
Lemma \ref{lem} it is enough to consider those greedy defining sets
of $G$ which are contained in $X$. Based on the property of our
ordering $\sigma$, we obtain that a descent in $G$ consists only of
a vertex from $E_2$ together with its two endpoints in $V_1$. This
shows that a greedy defining set of $G$ is a subset of $V_1$ which
dominates the elements of $E_2$, i.e. a vertex cover of $F$. The
converse is also true. It turns out that $GDN(G)$ is the same as the
minimum size of a vertex cover in $F$. This proves the theorem for
this case. The case where $X$ is colored by 2 is proved similarly in
which a subset $S\subseteq Y$ is a GDS for $G$ if and only if
$S\subseteq V_2$ and it dominates all elements of $E_1$. Namely in
this case too the minimum GDS is the same as the smallest vertex
cover of $F$. This completes the proof.
\end{proof}

\begin{thm}
There exists an efficient algorithm to determine the greedy defining
number of a forest.
\end{thm}

\noindent \begin{proof} It is enough to prove the theorem for trees,
since suppose that a forest $F$ consists of the connected components
$T_1, T_2, \ldots, T_k$. Then $GDN(G)= \sum_i GDN(T_i)$.

\noindent Now let $T$ be an ordered tree. It contains exactly two
proper colorings using two colors since it is connected. It is
enough to determine the greedy defining number of a 2-coloring of
$T$. Let a 2-coloring be given by a bipartition $(X,Y)$ of $V(T)$
where $X$ consists of vertices colored 1. Recall that a descent is
of the form a vertex colored 2 say $v$ together with its all
neighbors of $v$. These neighbors are colored 1 and have higher
order than $v$. Consider the subgraph $T'$ of $T$ induced by the
vertices of the descents in $T$. By Lemma \ref{lem} it is enough to
find a subset of vertices of color 1 with the minimum cardinality
which dominates all the vertices of color 2 in $T'$. Such a subset
will be denoted by $K$ and constructed gradually. Let $G$ be a
connected component of $T'$. Assume first that a vertex $v$ colored
2 has degree one in $G$ and let $u$ be its neighbor of color 1. Then
$\{v,u\}$ forms a descent and $u$ should be put in $K$. Delete now
$\{v,u\}\cup N$ from $G$, where $N$ is the neighbors of $u$. We do
the same for other similar vertices of $G$. After this stage we
obtain a subgraph of $G$ whose all leaves are colored 1. Note that
since any leaf of color 1 dominates only one vertex therefore it is
enough to consider only those vertices of color 1 which are not leaf
as possible elements to be put in $K$. Hence at this stage we remove
all leaves of color 1 from the graph. Since the earlier graph has no
any cycle then at each stage of our algorithm there is at least one
leaf. If it is colored by 2 then we put its neighbor (colored 1) in
our GDS $K$. Otherwise we remove it and continue until no vertex in
$G$ is left. We repeat this procedure for other connected components
of $T'$. This completes the proof.
\end{proof}

\section{Latin squares}

\noindent The minimum greedy defining number of any $n\times n$
Latin square is denoted by $g(n)$ in \cite{Z2,Z3} where it was shown
that $g(n)=0$ when $n$ is a power of two. The exact values of $g(n)$
for $n\leq 6$ were given in \cite{Z3} and for $n=7, 9, 10$ in
\cite{V}. But the complexity status of determining the greedy
defining number of Latin squares is still unknown. In the sequel we
present a method to obtain a greedy defining set in a Latin square.

\noindent For any $n\times n$ Latin square $L$ on $\{1,2, \ldots,
n\}$ we correspond three graphs $R(L)$, $C(L)$ and $E(L)$. Let $R$
be an arbitrary row of $L$. We first define a graph $G[R]$ on the
vertex set $\{1, \ldots, n\}$ as follows. Two vertices $i$ and $j$
with $1\leq i<j\leq n$ are adjacent in $G[R]$ if and only if (1) $j$
appears before $i$ in the row $R$ and (2) there is another entry $i$
in the same column of $j$ such that it comes after $j$ (i.e. lower
than $j$). In other words $i$ and $j$ are adjacent if and only if
they form a descent (jointly with an additional entry $i$). The
graph $R(L)$ is now defined the disjoint union $\cup G[R]$ on $n^2$
vertices where the union is taken over all $n$ rows of $L$. For any
column $C$ of $L$ we define $G[C]$ similarly. The graph $C(L)$
consists of the disjoint union of $G[C]$'s. Finally, in the sequel
we define $E(L)$. Let $e\in \{1, \ldots, n\}$ be any fixed entry.
There are $n$ entries equal to $e$ in $L$. First, a graph denoted by
$G[e]$ on these $n$ entries is defined in the following form. Two
entries $e_1$ and $e_2$ (which both are the same as $e$ but in
different rows and columns) are adjacent if and only if with an
additional entry they form a descent in $L$. The disjoint union of
$G[e]$'s form a graph which we denote by $E(L)$. The following
proposition is immediate.

\begin{prop}
A subset $D$ of entries of a Latin square $L$ is a GDS if $D$ is a
vertex cover for at least one of the graphs $R(L)$, $C(L)$ and
$E(L)$.\label{prop1}
\end{prop}

\noindent Proposition \ref{prop1} provides some upper bounds for the
greedy defining number of Latin squares. As an application, in the
following we present an upper bound for the greedy defining number
of any Latin square. We recall that according to Turan's theorem any
graph on $n$ vertices and with no clique of order $m$ has at most
$(m-2)n^2/(2m-2)$ edges.

\begin{thm}
Any $n\times n$ Latin square contains a GDS of size at most
$n^2-\frac{n\log 4n}{4}$.\label{prop2}
\end{thm}

\noindent \begin{proof} It is enough to find a vertex cover for
$E(L)$ of the desired cardinality. For this purpose consider the $n$
connected components of $E(L)$ which correspond to distinct entries
of $L$ i.e. $G[1], G[2], \ldots, G[n]$. We obtain an upper bound for
the vertex cover of each $G[i]$. The number of edges of $G[i]$ is
maximized when the $n$ entries of $i$ lie in the northeast-southwest
diagonal of $L$ and the maximum possible number of entries greater
than $i$ are placed in the top of this diagonal. It turns out that
in this case the graph has no more than $n(n-1)/2 - i(i-1)/2$ edges.
Assume that $G[i]$ has at most $f(i)$ independent vertices. Then the
complement of $G[i]$ which has at least $i(i-1)/2$ edges, does not
contain a clique of order $f(i)+1$. Using Turan's theorem we obtain
$$f(i) \geq \frac{n^2}{n^2-i(i-1)}.$$

\noindent If we write $i=n-j$ for some $0\leq j\leq n-1$ then
$$f(i) \geq \frac{n^2}{n(2j+1)-(j^2+j)}\geq \frac{n}{2j+1}.$$

\noindent This shows that $E(L)$ contains at least $n\sum_{0}^{n-1}
\frac{1}{2j+1}$ independent vertices. But from other side

$$2\sum_{j=0}^{n-1}\frac{1}{2j+1} \geq \sum_{k=1}^{2n+1} \frac{1}{k}
\geq 1+\frac{\log (2n+1)-1}{2}\geq \frac{\log 4n}{2}.$$

\noindent It turns out that $E(L)$ contains at least $(n\log 4n)/4$
independent vertices. Therefore it contains a vertex cover of no
more than $n^2-(n\log 4n)/4$ vertices. This completes the proof.
\end{proof}

\noindent Finally we mention an application of GDS of Latin squares
in secret sharing schemes. There is a known technique to design
secret sharing schemes using critical sets in Latin squares
\cite{CDS}. In a Latin square $L$, a subset of entries $S$ is said
to be a critical set if $L$ is the only Latin square which contains
$S$ as its partial Latin square. In the model presented in
\cite{CDS}, a dealer chooses a Latin square $L$ as the key of the
scheme and then obtains a collection of critical sets of $L$. Then
she shares the key among the participants in such a way that the
participants of any authorized set receive the entries of a critical
set $S$ and therefore by pooling their shares they can extend $S$ to
obtain the key uniquely. Since $S$ is a critical set then the Latin
square obtained, is nothing but $L$ itself. There are two serious
practical problems for this model. First, it is not easy to find an
arbitrary number of non-trivial critical sets in a random Latin
square. Second, once a set of participants pool their shares and
obtain a partial Latin square $S$, it is difficult to determine
whether it can be uniquely completed to a Latin square and then to
obtain the value of the key. Because it was known that, given a
partial Latin square $S$, whether $S$ has a unique completion to a
Latin square is an NP-complete problem \cite{CCS}. Now we propose
the same scenario as above but in stead of critical sets we use GDS.
The advantages of our model are as follow. First, there are known
and convenient techniques to obtain greedy defining sets in Latin
squares. In fact Theorem \ref{descent} as a general method and
Proposition \ref{prop1} as its refined version provide some tools to
construct arbitrary greedy defining sets in Latin squares (note that
to detect the descents of a square is an easy job). Second and more
importantly, when we have a partial Latin square in hand, it is very
easy using greedy coloring to find out whether it extends (uniquely)
to a given Latin square or not.

\section{Acknowledgment}

The author is thankful to John van Rees for pointing out a flaw in
the proof of Theorem 2 in \cite{Z1}.

\end{document}